# Risk and resampling under model uncertainty

Snigdhansu Chatterjee[1] and Nitai D. Mukhopadhyay[2]

*University of Minnesota and Virginia Commonwealth University*

**Abstract:** In statistical exercises where there are several candidate models, the traditional approach is to select one model using some data driven criterion and use that model for estimation, testing and other purposes, ignoring the variability of the model selection process. We discuss some problems associated with this approach. An alternative scheme is to use a model-averaged estimator, that is, a weighted average of estimators obtained under different models, as an estimator of a parameter. We show that the risk associated with a Bayesian model-averaged estimator is bounded as a function of the sample size, when parameter values are fixed. We establish conditions which ensure that a model-averaged estimator's distribution can be consistently approximated using the bootstrap. A new, data-adaptive, model averaging scheme is proposed that balances efficiency of estimation without compromising applicability of the bootstrap. This paper illustrates that certain desirable risk and resampling properties of model-averaged estimators are obtainable when parameters are fixed but unknown; this complements several studies on minimaxity and other properties of post-model-selected and model-averaged estimators, where parameters are allowed to vary.

## Contents

## 1. Introduction

In typical statistical applications, it is rare that a precise model is available to fit to the data. Selecting one model from several competing models is often the first step in the process. However, in the subsequent analysis, it is common to ignore the variability in the initial model selection. Two of the many consequences of ignoring modeling variability are (i) under-estimation of the variability of estimators and predictors, and (ii) erroneous inference and prediction, resulting from incorrectly computing the distributions of estimators and predictors. An alternative to selecting a model first and then computing an estimator under that model is to consider several models and appropriately average the estimators computed under these models.

---

[1] School of Statistics, University of Minnesota, 313 Ford Hall, 224 Church Street SE, Minneapolis, MN 55455, USA, e-mail: chatterjee@stat.umn.edu
[2] Department of Biostatistics, Virginia Commonwealth University, 730 E. Broad Street, Richmond, Virginia 23298, USA, e-mail: ndmukhopadhy@vcu.edu

Several studies have been published recently on the properties of post-model-selected and model-averaged estimators; see for example, [8], [23] and [24]. These studies are discouraging as they show that many nice properties associated with estimators under a known model vanish when there is model uncertainty. For example, Yang [23] shows that consistent model selection/averaging, and minimax-rate optimality cannot be simultaneously obtained. The review of Leeb and Pötscher [8] contains a discussion of several other problems with inference after model selection.

In view of these negative results, it seems desirable to scale down our expectations while working under model uncertainty, and strive for positive, if weaker, results. This may be achieved in one of two ways: we may either impose less stringent conditions on our estimators, or we may relax the criterion by which an estimator is evaluated. The latter is the goal of the present study.

The computation of an estimator is generally one of the early steps in a statistical exercise. Estimators of parameters are used for various purposes, notably for quantifying evidence for or against scientific hypotheses, obtaining interval estimates for the parameter under consideration, for prediction and forecasting, and for quantifying the accuracy of predictions and forecasts. These applications require knowledge about the distribution of the estimator, and knowledge about the risk associated with the usage of such estimators. In this paper, we concentrate on the risk behavior of a model-averaged estimator, and on approximating the distribution of a model-averaged estimator using the bootstrap.

In the first part of our study we show that under the traditional frequentist assumption that the parameters are fixed but unknown constants, the mean squared error in regression estimation under consistent model selection/averaging is bounded as a function of sample size. This complements Yang [23], where it was shown that a similar quantity cannot achieve minimax-rate optimality. Several of the negative results, including those of Yang [23], arise when a parameter is a known constant in a smaller model, while it is allowed to vary in a local neighborhood of that constant in a larger model. Recently, Hjort and Claeskens [5] studied model averaged estimators under a local parameter framework. Local parameters are ideal for mathematical development, but they are not reflective of statistical reality; see [17]. Indeed, as Hjort and Claeskens themselves remark in the rejoinder to the discussion of their paper, "a too literal belief in sample-size-dependent parameters would clash with Kolmogorov consistency and other requirements of natural statistical models." [5]. In view of this, it is meaningful to verify that estimators have reasonable risk behavior under consistent model selection/averaging when parameters are fixed constants. Our result also implies that *integrated risks* under consistent model selection/averaging are bounded, when integrals are taken with respect to any probability measure on the parameter space that does not depend on sample size.

In the second part of our study, in addition to the assumption that the parameters are fixed but unknown constants, we also weaken the consistency requirement of the model averaging procedure. In the terminology of Yang [23], a model selection/averaging scheme is *consistent* if it is asymptotically degenerate at the true model, when the true model is one of the candidate models. When the models are nested and several of them can correctly describe the data generation process, the most parsimonious correct model is taken as the true model. We call this *strong consistency*. We define a model selection/averaging scheme as *weakly consistent* if it selects or averages over all candidate models that correctly describe the data generation process. When only one model is correct, the strong and weak consistency requirements are identical; but if models are nested and several of them are correct, a weakly consistent scheme may distribute weights among all of them while a strongly consistent one is asymptotically degenerate at the smallest one. Recently, Leung and Barron [11] proposed a scheme of model averaging that results in nice risk behavior. Their scheme is an example of a weakly consistent procedure. We show that a particular choice of a weakly consistent model-averaged estimator has a distribution that can be approximated using the bootstrap.

In Section 2 we propose a simple linear regression model framework to study model uncertainty. We also discuss some of the properties of post-model-selection estimators that make them unsuitable for further applications, and also some properties of model-averaged estimators. This is followed in Section 3 with a discussion of mean squared error of the Bayesian model-averaged estimator. In Section 4 we propose a new adaptive, model-averaged estimator whose distribution may be consistently approximated using the bootstrap. A simulation example is discussed in Section 5. Finally, in Section 6 we discuss some aspects of our results, and point to some open issues relating to model uncertainty.

## 2. Issues with model selection or averaging

We select a simple regression framework for our study, which is the same as that used by [8], and similar to that of [24]. The observed data $\{(Y_t, \mathbf{x}_t = (x_{t1}, x_{t2})^T), t = 1, \ldots, n\}$, are modeled as

$$Y_t = \alpha x_{t1} + \beta x_{t2} + e_t, \tag{2.1}$$

where the $e_t$'s are independent, identically distributed $N(0, \sigma^2)$, $\sigma^2$ known. The design matrix $\mathbf{X}$ with rows given by $\mathbf{x}_t^T = (x_{t1}, x_{t2})$ is non-random. We denote the two columns of $\mathbf{X}$ as $X_1$ and $X_2$, the vector of errors as $\mathbf{e}$, and the vector of observations as $\mathbf{Y}$. The inner products and norms used below are the usual Euclidean ones. The notation $D$ is used for the determinant of the design matrix, thus $D = ||X_1||^2 ||X_2||^2 - <X_1, X_2>^2$. The unknown parameters in this model are $(\alpha, \beta)$. Model uncertainty surrounds the issue of whether or not $\beta = 0$. In this paper, for ease in presentation, we consider the problem of estimation of $\alpha$.

We make the standard assumption that $n^{-1}\mathbf{X}^T\mathbf{X} \to Q$ for a positive definite matrix $Q$. This, in particular, implies the standard design conditions

$$||X_1||^2 = O(n), \quad ||X_2||^2 = O(n), \tag{2.2}$$
$$<X_1, X_2> = O(n), \quad D = ||X_1||^2 ||X_2||^2 - <X_1, X_2>^2 = O(n^2). \tag{2.3}$$

We also assume that $n^{-1} <X_1, X_2> \not\to 0$ as $n \to \infty$, since without this restriction the effect of model uncertainty vanishes in this framework.

The true model, called $M_0$, may be described as

$$M_0 = \begin{cases} U \text{ (unrestricted)} & \text{if } \beta \neq 0; \\ R \text{ (restricted)} & \text{if } \beta = 0. \end{cases}$$

Under $U$, we adopt the ordinary least squares or maximum likelihood estimators $\widehat{(\alpha, \beta)} = (\mathbf{X}^T\mathbf{X})^{-1}\mathbf{X}^T\mathbf{Y}$. Our notation for these are $(\hat{\alpha}(U), \hat{\beta}(U))$. Under $R$, $\hat{\beta}(R) \equiv 0$, and the ordinary least squares or maximum likelihood estimator for $\alpha$ is $\hat{\alpha}(R) = [\sum_i x_{1i}^2]^{-1} \sum x_{1i} y_i$. Define $V_1 = \sigma^{-1} ||X_1||^{-1} <X_1, \mathbf{e}>$ and $V_2 = \sigma^{-1} D^{-1/2} ||X_1|| \{<X_2, \mathbf{e}> - ||X_1||^{-2} <X_1, X_2><X_1, \mathbf{e}>\}$, thus $\mathbf{V} = (V_1, \ V_2)^T \sim N(0, \mathbf{I}_2)$. In terms of $\mathbf{V}$, the estimators are

$$\begin{pmatrix} \hat{\alpha}(R) \\ \hat{\alpha}(U) \\ \hat{\beta}(U) \end{pmatrix} = \begin{bmatrix} \alpha + \beta ||X_1||^{-2} <X_1, X_2> + \sigma ||X_1||^{-1} V_1, \\ \alpha + \sigma ||X_1||^{-1} V_1 - \sigma ||X_1||^{-1} D^{-1/2} <X_1, X_2> V_2, \\ \beta + \sigma ||X_1|| D^{-1/2} V_2. \end{bmatrix}$$

The dichotomy between the *bias* of the restricted model $R$ and the *variance* of the unrestricted model $U$ can be clearly seen in the above formula. The restricted model estimator $\hat{\alpha}(R)$ has a bias factor $\beta ||X_1||^{-2} <X_1, X_2>$, which vanishes under $R$, while $\hat{\alpha}(U)$ has an extra factor of $\sigma ||X_1||^{-1} D^{-1/2} <X_1, X_2> V_2$ that inflates its variance relative to $\hat{\alpha}(R)$. Hence, model selection or model averaging is essentially a process of balancing bias and variance; see [20].

Let $\sigma_\beta$ be the standard deviation of $\hat{\beta}(U)$. This is a non-random, known number depending on $\sigma^2$ and $\mathbf{X}$. The following model selection criterion is used:

$$\hat{M} = \begin{cases} U & \text{if } |n^{-1/2} \sigma_\beta^{-1} \hat{\beta}(U)| > c; \\ R & \text{if } |n^{-1/2} \sigma_\beta^{-1} \hat{\beta}(U)| \leq c. \end{cases}$$

The above criterion may be identified as representative of standard model selection tools, in the simple regression model. In particular, the above criterion is the traditional pre-test procedure based on the likelihood ratio, coincides with the *Akaike Information Criterion* (AIC) if $c = \sqrt{2}$, and coincides with the *Bayesian Information Criterion* (BIC) if $c = \sqrt{\log n}$. The post-model-selection estimator of $\alpha$ is

(2.4) $$\tilde{\alpha} = \hat{\alpha}(R) I_{\{\hat{M}=R\}} + \hat{\alpha}(U) I_{\{\hat{M}=U\}}.$$

Several nice properties are known about $\hat{M}$ and, consequently, it is generally believed that $\tilde{\alpha}$ will also have good properties. Some of the important properties include that for all $\beta$ and as $c \to \infty, n^{-1/2} c \to 0$, $P[\hat{M} = M_0] \to 1$, $\{\hat{M} = M_0\} \subseteq \{\tilde{\alpha} = \hat{\alpha}(M_0)\}$ and thus $P[\tilde{\alpha} = \hat{\alpha}(M_0)] \to 1$ (see [15]). Note that $\hat{\alpha}(M_0)$ is the "oracle's guess" about $\alpha$, and is not a statistic, since it is based on the knowledge of $\beta$. The above properties tend to give the impression that $\tilde{\alpha}$ is a very good estimator.

However, there are some major problems since the above results are asymptotic in nature, and the asymptotics can take a long time to kick in, as well as be dependent on the value of $\beta$. Our primary reference for this model and its basic properties [8] identifies this as a problem of *non-uniformity* in $\beta$ of the convergence of $\hat{M}$ and $\tilde{\alpha}$. It can be immediately seen that the estimator $\tilde{\alpha}$ is super-efficient when $c \to \infty$, $c/\sqrt{n} \to 0$, as with BIC. The major repercussions of super-efficiency of $\tilde{\alpha}$ and the non-uniformity of its asymptotics is in its risk performance, and in its finite sample behavior. The mean squared error of $\tilde{\alpha}$ is unbounded and depends on $\beta$, while that of $\hat{\alpha}(M_0)$ is a constant. As a consequence, the finite sample behavior of $\tilde{\alpha}$ is erratic and can be quite unlike its asymptotic approximation. Available simulations confirm this; see [8]. Several other studies conducted by Leeb, Pötscher, Yang and others reveal how and why the properties of $\tilde{\alpha}$ and $\hat{\alpha}(M_0)$ differ. For further information see, for example, [6, 7, 8, 9, 10, 22, 24, 25].

The super-efficiency of $\tilde{\alpha}$ results in most variations of the bootstrap being inapplicable. Only subsampling ([14]) and the $m$-out-of-$n$ bootstrap with $m/n \to 0$ would yield consistent approximations of the distribution of $\tilde{\alpha}$. Unfortunately, these methods have problems of their own, some details of which can be found in [18] and [1]. Specifically, although subsampling is asymptotically consistent, it can perform miserably in finite samples. For any $\alpha \in (0, 1)$, the actual asymptotic coverage of a standard level $(1 - \alpha)$ subsampling confidence interval can be zero; see [1] for details. The finite sample properties of subsampling based methods can be improved sometimes by considering hybrid techniques, calibrations and other modifications, as documented by [2]. However, the asymptotic zero coverage of subsampling intervals for $\tilde{\alpha}$ cannot be reversed by, for example, size correction, since technical conditions that allow for such correction to work are not satisfied by $\tilde{\alpha}$.

The above issues with post-model-selection estimators lead to model-averaged estimators. A model-averaged estimator of $\alpha$ is of the form

$$\check{\alpha} = \hat{\alpha}(R) p_R + \hat{\alpha}(U) p_U, \tag{2.5}$$

where $p_R$ and $p_U$ are two *weights* associated with the models $R$ and $U$. Yang and his co-authors have extensively studied aggregation across models for several statistical procedures like estimators and forecasts, in both their algorithmic as well as theoretical aspects (see [22, 23, 24, 25]). In particular, a result of [23] implies that when the model averaging technique is strongly consistent, the supremum of the mean squared error of $n^{1/2}(\check{\alpha} - \alpha)$ over values of $(\alpha, \beta)$ tends to infinity. Thus, strongly consistent model averaging does not attain the minimax rate. Our result in Section 3 shows that, up to constant terms, it is no worse than the post-model-selection estimator when $(\alpha, \beta)$ are held fixed.

Recently, [5] studied several forms of model averaging and showed that a typical model-averaged estimator converges weakly to a mixture of normal laws, when the parameters of the true model are in a $O(n^{-1/2})$ neighborhood of the simplest candidate in a nesting of models. Since subsampling does not seem to perform well in practice, it is important to study conditions on model weights under which bootstrap approximations of finite sample distributions hold, i.e., conditions under which the statistic under consideration is smooth and asymptotically normal (see [12], [13]). This is studied in Section 4.

## 3. Risk profile of model-averaged estimators

Several problems associated with the post-model-selection estimator can be attributed to its lack of uniformity, as discussed extensively by others [8]. One is the super-efficiency of $\tilde{\alpha}$, for example, when BIC is used for model selection. The core problem of lack of uniformity in the convergence pattern of $\tilde{\alpha}$ is unavoidable – even with model averaging – when a strongly consistent model averaging technique is used, as described by [23]. In this section we show that when parameter values are fixed, model averaging is no worse than model selection, up to constant terms.

Under the unrestricted model, $U$, we choose the prior on $(\alpha, \beta)$ to be a standard mean zero, identity covariance bivariate Normal distribution, $N(0, I)$. Under the restricted model, $R$, the prior on $\alpha$ is a standard univariate Normal distribution, $N(0, 1)$. We put equal prior weights, i.e., $1/2$, on the models, so the prior odds is 1. Our notation for the posterior probabilities of the two models are $\pi_{nU}$ and $\pi_{nR}$. Since $\sigma$ is known, without loss of generality we also assume $\sigma = 1$ in this section.

Thus the Bayesian model-averaged estimator of $\alpha$ is

(3.6) $$\hat{\alpha}_{BMA} = \pi_{nU}\hat{\alpha}(U) + \pi_{nR}\hat{\alpha}(R).$$

We use the pre-selected, least squares estimators $\hat{\alpha}(U)$ and $\hat{\alpha}(R)$ as constituents of $\hat{\alpha}_{BMA}$, and consider the squared error loss function. The case where a general loss function is used, with $\hat{\alpha}(U)$ and $\hat{\alpha}(R)$ taken to be the Bayes estimators under models $U$ and $R$, is very similar. The following Proposition is our main result in this section.

**Proposition 3.1.** *The normalized risk of $\hat{\alpha}_{BMA}$, $nR(\alpha) = nE(\hat{\alpha}_{BMA} - \alpha)^2$, satisfies $\sup_n nR(\alpha) < \infty$, for every fixed choice of $\alpha$ and $\beta$. Hence, the integrated normalized risk*

$$\sup_n \int_{\alpha,\ \beta} nR(\alpha) d\lambda(\alpha,\ \beta) < \infty$$

*for any probability measure $\lambda(\cdot)$ that does not depend on $n$.*

*Proof.* In the following, we use $C$ as a generic constant, not depending on the parameters $\alpha$ and $\beta$ or the sample size $n$.

Note that $\hat{\alpha}(R) = \hat{\alpha}(U) + \hat{\beta}(U)||X_1||^{-2} <X_1, X_2>$. Therefore,

$$nR(\alpha) = nE\left[\pi_{nU}\hat{\alpha}(U) + \pi_{nR}\hat{\alpha}(R) - \alpha\right]^2$$
(3.7) $$\leq 2nE(\hat{\alpha}(U) - \alpha)^2 + 2n||X_1||^{-4} <X_1, X_2>^2 E\left\{\pi_{nR}^2\hat{\beta}^2(U)\right\}.$$

Note that $E(\hat{\alpha}(U) - \alpha)^2 = \sigma^2||X_1||^{-2}E\left[V_1 - <X_1, X_2> D^{-1/2}V_2\right]^2 = Cn^{-1}$ and $E\pi_{nR}^2\hat{\beta}^2(U). \leq 2\beta^2 E\pi_{nR}^2 + Cn^{-1}$. Thus, we need suitable bounds for $\beta^2 E\pi_{nR}^2$. We now have

$$p_{nR} = m_R(\mathbf{Y})/(m_U(\mathbf{Y}) + m_U(\mathbf{Y})) = \frac{m_R(\mathbf{Y})}{m_U(\mathbf{Y})}\left(1 + \frac{m_R(\mathbf{Y})}{m_U(\mathbf{Y})}\right)^{-1} \leq \frac{m_R(\mathbf{Y})}{m_U(\mathbf{Y})}.$$

Then, making use of the moment generating function of a $\chi^2$ random variable, we can deduce that

$$E\left(\frac{m_R(\mathbf{Y})}{m_U(\mathbf{Y})}\right)^2 = Cn^2 \exp\left\{-nC_0(\alpha^2 + \beta^2)\right\}$$

for a particular constant $C_0$. This yields, at (3.7), that

$$nR(\alpha) = Cn^{-1} + Cn^3\beta^2 \exp\left\{-nC_0(\alpha^2 + \beta^2)\right\}.$$

which is bounded for every fixed $(\alpha,\ \beta)$, as a function of $n$. The rest of the result follows. $\square$

**Remark 3.1.** A lower bound for $nR(\alpha)$ can also be established using arguments similar to those above. With slight modification, the above approach using the moment generating function of a non-central $\chi^2$ random variable can be used to provide an alternative proof of Theorem 2 of [23]. It can also be seen that even when $(\alpha,\ \beta)$ vary over a compact set, the supremum of $nR(\alpha)$ over $(\alpha, \beta)$ is unbounded.

## 4. Adaptive model-averaged estimators and the bootstrap

The results of Hjort and Claeskens [5] and Leeb and Potscher [8] indicate that the post-model-selection estimator and many model-averaged estimators cannot be consistently bootstrapped. The problems associated with the risk behavior, and those associated with bootstrap approximation, arise from two different sources. Undesirable behavior of the risk function arises from considering scenarios as parameters vary, while a major reason why the distribution of post-model-selection or model-averaged estimators cannot be approximated by bootstrap methods is because of lack of smoothness of the estimator, or lack of asymptotic normality.

In this section we study the conditions on the model weights which are required for consistent bootstrap approximation of the distribution of the resulting model-averaged estimator. Clearly, since the distribution of $\hat{\alpha}(U)$ can be approximated using the bootstrap, putting the entire weight on model $U$ is an option. However, balancing between $\hat{\alpha}(U)$ and $\hat{\alpha}(R)$ can lead to a more efficient estimator. We propose below a data-adaptive model weighing scheme that achieves the dual goals of reasonable efficiency and bootstrap consistency.

A model-averaged estimator of $\alpha$ is of the form

$$\check{\alpha} = \hat{\alpha}(R) p_{nR} + \hat{\alpha}(U) p_{nU}. \tag{4.8}$$

Notice that we have adopted a different notation ($p_{nR}$ and $p_{nU}$) for the model weights in this Section, from those ($\pi_{nR}$ and $\pi_{nU}$) used in Section 3. This is to emphasize that the nature of these weights may be different. We retain the condition that the parameters $(\alpha, \beta)$ are fixed but unknown.

A primary requirement for consistency is $p_{nR} + p_{nU} = 1$, as pointed out in [5]. In order to avoid pathologies, we also specify that $p_{nU} \in [0, 1]$. Note that the weights $p_{nR}$ and $p_{nU}$ may depend on the parameters $(\alpha, \beta)$, and the random component $\mathbf{V}$, apart from the known constants $\mathbf{X}$ and $\sigma^2$.

Replacing $p_{nU}$ by $1 - p_{nR}$, we thus have

$$\begin{aligned}\check{\alpha} &= \alpha + \sigma ||X_1||^{-1} V_1 + \beta p_{nR} ||X_1||^{-2} <X_1, X_2> \\ &\quad - \sigma ||X_1||^{-1} D^{-1/2} <X_1, X_2> (1 - p_{nR}) V_2.\end{aligned}$$

A primary requirement on $\check{\alpha}$ is that it should be consistent, and the following proposition establishes a necessary and sufficient condition for this.

**Proposition 4.1.** *The model-averaged estimator $\check{\alpha}$ converges in probability to $\alpha$ if and only if $\beta p_{nR}$ converges in probability to zero as $n \to \infty$.*

*Proof.* The sufficiency part follows easily from the design conditions (2.2)–(2.3). For the necessity part, suppose that $\beta p_{nR} \xrightarrow{p} \tilde{c} \neq 0$ as $n \to \infty$. This is clearly equivalent to $p_{nR} \xrightarrow{p} c = \tilde{c}/\beta \neq 0$ as $n \to \infty$ and $\beta \neq 0$. Hence, we also have $(1 - p_{nR}) \left\{ \sigma ||X_1||^{-1} D^{-1/2} <X_1, X_2> V_2 \right\} \xrightarrow{p} (1-c) 0 = 0$. This implies $\check{\alpha} \xrightarrow{p} \alpha - \tilde{c}\gamma \neq \alpha$, where $||X_1||^{-1} <X_1, X_2> \to \gamma$ as $n \to \infty$. The case where $p_{nR}$ does not have a limit can be treated similarly with a little more algebra. $\square$

The next proposition is an extension of the previous one, and establishes sufficient conditions for asymptotic normality of $\check{\alpha}$.

**Proposition 4.2.** *The scaled and centered model-averaged estimator $n^{1/2}(\check{\alpha} - \alpha)$ has an asymptotic normal distribution if* (i) *$n^{1/2} \beta p_{nR}$ converges in probability to zero as $n \to \infty$, and* (ii) *$p_{nR}$ converges in probability as $n \to \infty$ for all values of $(\alpha, \beta)$.*

*Proof.* The first condition forces the bias component in $\check{\alpha}$ to be $o(n^{-1/2})$, while the second condition allows for use of Slutsky's theorem. □

By requiring $n^{1/2}\beta p_{nR} \xrightarrow{p} 0$ as $n \to \infty$ we have ensured that, when $\beta \neq 0$, we have $n^{1/2}p_{nR} \xrightarrow{p} 0$. Thus the model-averaged estimator is close to the unrestricted model estimator $\hat{\alpha}(U)$, and has the same limiting distribution up to first order terms. However, when $\beta = 0$, the asymptotic distribution of $n^{1/2}(\check{\alpha} - \alpha)$ depends on the limit of $p_{nR}$, which is between zero and one. Thus, when the restricted model holds, the asymptotic variance of $\check{\alpha}$ is between that of $\hat{\alpha}(R)$ and $\hat{\alpha}(U)$. The relative strengths of different candidates for model weight $p_{nR}$ may be evaluated by their probability limits when $\beta = 0$. We note that we consider $(\alpha, \beta)$ as fixed constants and do not allow them to vary with $n$. If, for example, we assumed $\beta = O(n^{-1/2})$, then the first condition of Proposition 4.2 would imply asymptotically zero weight on the restricted model.

In order to progress towards bootstrap consistency, apart from asymptotic normality of $\check{\alpha}$, we also need $p_{nR}$ to be a smooth function. Thus ruling out the indicator function $p_{nR} = I_{\{|n^{-1/2}\sigma_\beta^{-1}\hat{\beta}(U)| \leq c\}}$ used in $\tilde{\alpha}$. Keeping in view the nice properties of $\tilde{\alpha}$, we now develop an adaptive, data-driven model weight function $p_{nR}$ that is a smooth version of $I_{\{|n^{-1/2}\sigma_\beta^{-1}\hat{\beta}(U)| \leq c\}}$.

For any $k_n$, we split the event $\{-k_n \leq \hat{\beta}(U) \leq k_n\}$ into two events, $\{\hat{\beta}(U) - k_n \leq 0\}$ and $\{\hat{\beta}(U) + k_n \geq 0\}$, and approximate the indicators of these events separately. Our approximation for $I_{\{\hat{\beta}(U) - k_n \leq 0\}}$ is

$$\begin{aligned}\xi_{1n} &\equiv \xi_{1n}\left(\gamma_{1n}, \hat{\beta}(U), k_n\right) \\ &= \left(1 + \exp\left\{-\gamma_{1n}(\hat{\beta}(U) - k_n)\right\}\right)^{-1} \exp\left\{-\gamma_{1n}(\hat{\beta}(U) - k_n)\right\},\end{aligned}$$

and for $I_{\{\hat{\beta}(U) + k_n \geq 0\}}$ is

$$\begin{aligned}\xi_{2n} &\equiv \xi_{2n}\left(\gamma_{1n}, \hat{\beta}(U), k_n\right) \\ &= \left(1 + \exp\left\{\gamma_{2n}(\hat{\beta}(U) + k_n)\right\}\right)^{-1} \exp\left\{\gamma_{2n}(\hat{\beta}(U) + k_n)\right\}.\end{aligned}$$

We take the two tuning values $\gamma_{1n}$ and $\gamma_{2n}$ to be always positive. However, they change with $n$; and in a major departure from traditional model weights, they are not equal to each other, and also depend on the data. Thus, $\gamma_{1n} \equiv \gamma_{1n}(\alpha, \beta, \mathbf{V})$ and $\gamma_{2n} \equiv \gamma_{2n}(\alpha, \beta, \mathbf{V})$ are unequal, random weights.

Equipped with these functions, we define

$$p_{nR} = 0.5\xi_{1n} + 0.5\xi_{2n}.$$

We adopt the *paired bootstrap* as our resampling strategy. Thus, we draw a simple random sample with replacement of the data pairs $(Y_i^*, \mathbf{x}_i^*)$, $i = 1, \ldots, n$, from the original data $(Y_i, \mathbf{x}_i)$, $i = 1, \ldots, n$. The entire process of obtaining $\hat{\alpha}(R)$, $\hat{\alpha}(U)$, $\hat{\beta}(R)$, $p_{nR}$, and $\check{\alpha}$ is imitated with the resample $(Y_i^*, \mathbf{x}_i^*)$, $i = 1, \ldots, n$, and we approximate the distribution of $n^{1/2}(\check{\alpha} - \alpha)$ with the distribution of $n^{1/2}(\check{\alpha}^* - \check{\alpha})$, conditional on $(Y_i, \mathbf{x}_i)$, $i = 1, \ldots, n$. A technical condition guarantees that the design matrix from the resampled data is non-singular with high probability; see condition (1.17) of [3].

The following Theorem is our main result in this section, and establishes consistency of the bootstrap for a adaptively weighted model-averaged estimator.

**Theorem 4.1.** *Assume that sequence of constants $k_n \downarrow 0$ as $n \to \infty$. Suppose the tuning constants are chosen as $\gamma_{1n} = a_n \hat{\beta}(U)$ $\gamma_{2n} = -a_n \hat{\beta}(U)$ where $\{a_n\}$ is a sequence of positive constants satisfying $a_n^{-1} \log(n) \downarrow 0$ as $n \to \infty$.*

*Then $n^{1/2}(\check{\alpha} - \alpha)$ has an asymptotic Normal distribution and the paired bootstrap is consistent for it.*

*Proof.* For the asymptotic normality we only need to check that the conditions of Proposition 4.2 are met. We illustrate the calculation for verifying $n^{1/2} \xi_{1n} \xrightarrow{p} 0$, when $\beta \neq 0$.

$$P\left[|n^{1/2}\xi_{1n}| > \epsilon\right]$$
$$= P\left[|\left(1 + \exp\left\{-\gamma_{1n}(\hat{\beta}(U) - k_n)\right\}\right)^{-1}\right.$$
$$\left.\exp\left\{-\gamma_{1n}(\hat{\beta}(U) - k_n) + 0.5\log(n)\right\}| > \epsilon\right]$$
$$\leq P\left[\exp\left\{-\gamma_{1n}(\hat{\beta}(U) - k_n) + 0.5\log(n)\right\} > \epsilon\right]$$
$$= P\left[\hat{\beta}(U) \text{ lies between the roots of } x^2 - k_n x - 0.5 a_n^{-1} \log(n) + a_n^{-1} \log(\epsilon) = 0\right].$$

The roots of the equation $x^2 - k_n x - 0.5 a_n^{-1} \log(n) + a_n^{-1} \log(\epsilon) = 0$ are always real when $\epsilon < 1$, since $k_n^2 + 2a_n^{-1} \log(n) - 4a_n^{-1} \log(\epsilon) > 0$ for all $n$. Note that the square of the distance between the roots is given by $\left(k_n^2 + 2a_n^{-1} \log(n) - 4a_n^{-1} \log(\epsilon)\right)/4$. When $k_n \downarrow 0$, $k_n^2 + 2a_n^{-1} \log(n) - 4a_n^{-1} \log(\epsilon) \downarrow 0$, hence the Lebesgue measure of the interval between the roots goes to zero as $n \to \infty$, thus ensuring

$$P\left[\hat{\beta}(U) \text{ lies between the roots of } x^2 - k_n x - 0.5 a_n^{-1} \log(n) + a_n^{-1} \log(\epsilon) = 0\right] \to 0,$$

as $n \to \infty$. Note that this result actually does not depend on the value of $\beta$, as long as it is non-zero.

Other parts of the proof for asymptotic Normality may be verified similarly. Since $\check{\alpha}$ is a smooth function of $\alpha$, $\beta$ and $\mathbf{V}$, and has an asymptotic Normal distribution, the consistency of the paired bootstrap procedure follows from [12] and [13]. □

**Remark 4.1.** The condition $k_n \downarrow 0$ as $n \to \infty$ is a weaker restriction than typically found in literature. Since $\hat{\beta}(U) = O_p(n^{-1/2})$, the AIC criterion uses $k_n = O(n^{-1/2})$, while the BIC uses $k_n = O(n^{-1/2}\sqrt{\log(n)})$.

**Remark 4.2.** The assumptions of Proposition 4.1 and Proposition 4.2 cannot be weakened in general. The example of Section 10.6 of [5] provides a test case. It is a simpler version of the model described in Section 2, and simply has $Y_1, \ldots, Y_n$ independent, identically distributed as $N(\mu, 1)$ random variables. Model uncertainty is about whether $\mu = 0$, and the natural estimator for $\mu$ is $\bar{Y}_n = n^{-1} \sum_{i=1}^n Y_i$ in the unrestricted model, and 0 in the restricted model. A model-averaged estimator is $\hat{\mu} = W(n^{1/2}\bar{Y}_n)\bar{Y}_n$, for some weight $W(\cdot) \in [0,1]$. Note that under a model with contiguous alternatives $\mu_{\text{true}} = n^{-1/2}\delta$, the requirement that $\hat{\mu}$ be consistent for $\mu_{\text{true}}$ actually places no restriction on the weight $W(\cdot)$, which may take any value in $[0,1]$. However, if we want consistency under arbitrary $\mu$, $W(n^{1/2}\bar{Y}_n) \xrightarrow{p} 1$ is a requirement.

For asymptotic normality, $n^{1/2}\mu(1 - W(n^{1/2}\bar{Y}_n)) \xrightarrow{p} 0$ and convergence in probability of $W(n^{1/2}\bar{Y}_n)$, are requirements. Under $\mu_{\text{true}}$, this implies that $W(n^{1/2}\bar{Y}_n) \xrightarrow{p}$ 1 must hold, while for general $\mu$, the stronger condition $n^{1/2}(1 - W(n^{1/2}\bar{Y}_n)) \xrightarrow{p} 0$ must be satisfied.

Under $\mu_{\text{true}}$, it is of interest to approximate the distribution of the standardized statistic

$$\Lambda_n = n^{1/2}(\hat{\mu} - \mu_{\text{true}}) = n^{1/2}W(n^{1/2}\bar{Y}_n)\bar{Y}_n - \delta = W(\delta + Z_n)(\delta + Z_n) - \delta,$$

where $Z_n \sim N(0,1)$.

A natural question is what should be a bootstrap equivalent of $\Lambda_n$. Suppose $Y_1^*, \ldots, Y_n^*$ are a random sample from the data $Y_1, \ldots, Y_n$. We consider the bootstrap equivalent of $n^{1/2}\bar{Y}_n$ to be $n^{1/2}(\bar{Y}_n^* - \bar{Y}_n)$, and not $n^{1/2}\bar{Y}_n^*$. This is in keeping with [4], who put forth the guideline that for good power performance, resampling must be done to reflect the null hypothesis. While model selection is not in general a hypothesis test, some of the same principles are applicable.

Hence, we have $\hat{\mu}^* = W(n^{1/2}(\bar{Y}_n^* - \bar{Y}_n))\bar{Y}_n^*$. When $1 - W(n^{1/2}\bar{Y}_n) \xrightarrow{p} 0$, it can be readily seen that the distribution of $\Lambda_n^* = n^{1/2}(\hat{\mu}^* - \hat{\mu})$, conditional on $Y_1, \ldots, Y_n$, and that of $\Lambda_n$ converge to the same limit law. □

**Remark 4.3.** We conjecture that for the model-averaged estimator proposed in this section, a result similar to [16] would hold. In the framework of this paper, the statement corresponding to the main result of [16] would be as follows: Let $F_{n,\alpha,\beta}(t) = P\left[n^{1/2}(\check{\alpha} - \alpha) \leq t\right]$, and let $\hat{F}_n(t)$ be an estimator of $F_{n,\alpha,\beta}(t)$ satisfying for every $\delta > 0$ $P_{n,\alpha,\beta}[|\hat{F}_n(t) - F_{n,\alpha,\beta}(t)| > \delta] \to 0$, as $n \to \infty$. Then $\exists \delta_0 > 0$ and $\rho_0 > 0$ such that

$$(4.9) \quad \sup_{(\tilde{\alpha},\tilde{\beta}) \in B((\alpha,\beta);\rho_0/\sqrt{n})} P_{n,\tilde{\alpha},\tilde{\beta}}\left[|\hat{F}_n(t) - F_{n,\tilde{\alpha},\tilde{\beta}}(t)| > \delta_0\right] \to 1;$$

$$\text{where} \quad B((\alpha,\beta);a) = \{(\tilde{\alpha},\tilde{\beta}) : ||(\tilde{\alpha},\tilde{\beta}) - (\alpha,\beta)|| < a$$

is the open ball of radius $a$ around $(\alpha, \beta)$. It can be seen that under standard conditions, if the supremum in (4.9) is taken over $B((\alpha, \beta); a_n)$ with $a_n = o(n^{-1/2})$ instead of $B((\alpha, \beta); \rho_0/\sqrt{n})$, the limit would be zero instead of 1. Thus the result of [16] may be improved to the case where the supremum is taken only over the set of parameter values that are exact order $n^{-1/2}$ away from the $(\alpha, \beta)$ under which the estimator $\hat{F}_n(\cdot)$ is computed. This is easily verified, for example, when $\alpha = 0$, $\sigma = 1$ and $X_{t2} \equiv 1$.

Note that from a bootstrap approximation point of view, (4.9) is not a negative result, but a very positive one. The uses of bootstrap approximation are for constructing interval estimates, testing hypotheses and so on. Equation (4.9) and other related results from [16] imply that a bootstrap approximation $\hat{F}_n(\cdot)$ constructed under the "null" $(\alpha, \beta)$, has sup-norm distance of 1 from the true distributions under parameter values that are exact order $n^{-1/2}$ away from the $(\alpha, \beta)$. Thus $\hat{F}_n(\cdot)$ has power 1 in hypothesis testing under contiguous alternatives. This is a further confirmation of the tenet of [4], that resampling procedure ought to reflect the null hypothesis.

**Remark 4.4.** It is of interest to know that the asymptotic variance of $\check{\alpha}$ depends on $\beta$, and is given by $\text{Var}(n^{1/2}(\check{\alpha} - \alpha)) - \text{Var}(n^{1/2}(\hat{\alpha}(U) - \alpha)) \to 0$ if $\beta \neq 0$, while $\text{Var}(n^{1/2}(\check{\alpha}-\alpha)) - \{0.5\,\text{Var}(n^{1/2}(\hat{\alpha}(U)-\alpha)) + 0.5\,\text{Var}(n^{1/2}(\hat{\alpha}(R)-\alpha))\} \to 0$ if $\beta = 0$. This is established by checking that both $\xi_{1n}$ and $\xi_{2n}$ tend to $1/2$ as $n \to \infty$ when $\beta = 0$. Thus $\check{\alpha}$ performs like the correct estimator $\hat{\alpha}(U)$ when model $U$ is valid, and balances between the correct and conservative choices when the restricted model $R$ is true.

## 5. A simulation example

We performed a small simulation experiment to illustrate some of the features of inference under model uncertainty that have been discussed in the previous sections. We took $n = 50$, $x_{i1} \equiv 1$, and generated 50 numbers from the Uniform distribution supported between zero and three and fixed these as the $x_{i2}$ values. We fixed $\alpha = 1$, and varied the $\beta$ values.

For different values of $\beta \in [-1, 1]$, we obtained sampling distribution approximations of (i) the post-model-selected estimator $\hat{\alpha}_{MS}$, (ii) a version of the Bayesian model-averaged estimator $\hat{\alpha}_{BMA}$, and (iii) an adaptive model-averaged estimator $\hat{\alpha}_{AMA}$, by 5000 replications for each value of $\beta$. For the Bayesian model-averaged estimator, model $R$ was assigned weight $q_{nR} = \exp(-BIC_R/2)/(\exp(-BIC_R/2) + \exp(-BIC_U/2))$ while model $U$ was assigned weight $1 - q_{nR}$. We define

$$BIC_R = \sum [Y_i - \hat{\alpha}_R x_{i1}]^2 + \log(n),$$
$$BIC_U = \sum \left[Y_i - \hat{\alpha}_U x_{i1} - \hat{\beta}_U x_{i1}\right]^2 + 2\log(n).$$

For the adaptive model-averaged estimator, we took $a_n = (\log(n))^2$.

The requirement that $a_n^{-1} \log(n) \downarrow 0$ suggests that $a_n$ should be an increasing sequence, growing faster than $\log(n)$. Several choices of $a_n$ were used initially, and it turned out that very slowly increasing sequences like $a_n = (\log(n))^2$ or very quickly increasing sequences like $a_n = n^{0.499}$ performed better than others. This is a reflection on our way of constructing the functions $\xi_{1n}$ and $\xi_{2n}$ using $\gamma_{1n}$ and $\gamma_{2n}$. Alternative choices, like $\gamma_{1n} = a_n |\hat{\beta}(U)| \{\hat{\beta}(U)\}^{-1}$, are a subject for further research.

The first object of our study is the mean squared error of the three estimators of $\alpha$, namely, $\hat{\alpha}_{MS}$, $\hat{\alpha}_{BMA}$, and $\hat{\alpha}_{AMA}$. Panel (a) in Figure 1 contains the graphs of the mean squared error (MSE) as $\beta$ varies between $[-1, 1]$. In this and all subsequent figures, the solid line corresponds to $\hat{\alpha}_{BMA}$, the broken line to $\hat{\alpha}_{MS}$, and the dotted line to $\hat{\alpha}_{AMA}$. In this figure, we have also added the graph for the MSE of $\hat{\alpha}(U)$, which is the nearly horizontal dot-and-dash line. First, using model selection or averaging is clearly better than using $\hat{\alpha}(U)$ only in the region $0 \pm 2/\sqrt{n} \approx (-0.3, 0.3)$, where $MS$, $BMA$ and $AMA$ all perform better than $\hat{\alpha}(U)$. However, in the neighboring regions $|\beta| \in (0.3, 0.8)$, $\hat{\alpha}(U)$ has smaller MSE than the three estimators. For high values of $|\beta|$, using model selection/averaging or the unrestricted model makes little difference. Thus whether model averaging/selection is useful or not depends considerably on the value of $\beta$. Also note that $BMA$ has a lower MSE compared to $MS$ for low values of $|\beta|$ and only marginally higher MSE otherwise, with a much lower maximum MSE value. The graph for $AMA$ tends to stay closest to the graph for $\hat{\alpha}(U)$, and thus does better than $BMA$ or $MS$ in the region $|\beta| \in (0.05, 0.75)$, but is marginally poorer otherwise.

In order to study how the three estimators balance between $\hat{\alpha}(R)$ and $\hat{\alpha}(U)$, we computed the Kolmogorov–Smirnov distances $KS_{jR}$ and $KS_{jU}$, between the distribution of $n^{1/2}(\hat{\alpha}_j - \alpha)$, and the distributions of $n^{1/2}(\hat{\alpha}_R - \alpha)$ and $n^{1/2}(\hat{\alpha}_U - \alpha)$, where $j = MS, BMA, AMA$ (MS: model selected, BMA=Bayesian model-averaged, AMA=adaptive model-averaged). We then computed the ratios

$$KSRatio_j = 100 \frac{KS_{jR}}{KS_{jR} + KS_{jU}}, \quad j = MS, BMA, AMA.$$

Under ideal circumstances, this ratio ought to be zero at $\beta = 0$, and 100 for $\beta \neq 0$. Panel (b) in Figure 1 displays the $KSRatio_j$ values for the three estimators $j =$ $MS, BMA, AMA$. When $\beta = 0$, $MS$ is closest to $\hat{\alpha}_R$, while, as predicted, $AMA$ balances between $\hat{\alpha}_R$ and $\hat{\alpha}_U$. The Bayesian model-averaged estimator $BMA$ lies between $MA$ and $AMA$, and is quite close to $MS$. In the region $0 \pm 2/\sqrt{n} \approx (-0.3, 0.3)$ both $MS$ and $BMA$ are much closer to $\hat{\alpha}_R$ than $\hat{\alpha}_U$.

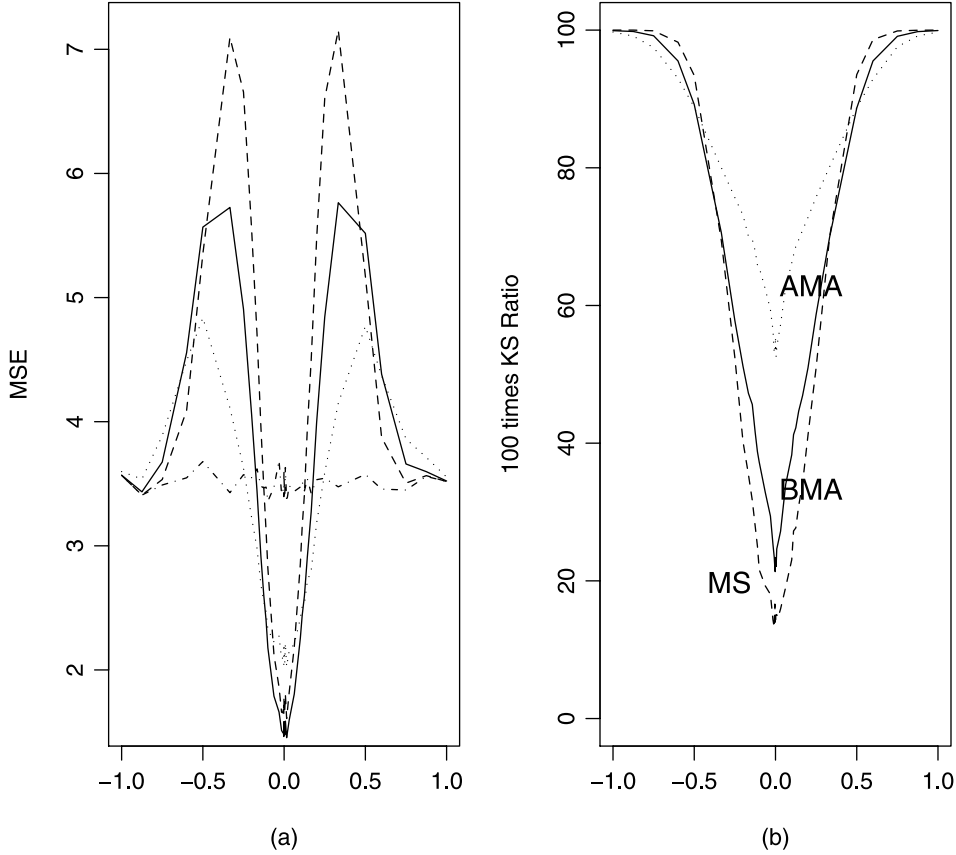

FIG 1. *Panel* (a) *is the mean squared error of $\hat{\alpha}_{BMA}$ (solid line), $\hat{\alpha}_{MS}$ (broken line), $\hat{\alpha}_{AMA}$ (dotted line), and $\hat{\alpha}(U)$ (dot-and-dash line). Panel* (b) *is the ratio of Kolmogorov Smirnov distances $KSRatio_j = KS_{jR}/(KS_{jR} + KS_{jU})$, $j = MS, BMA, AMA$, scaled by 100; between distributions of centered and scaled estimators and $\hat{\alpha}(R)$ (for $KS_{jR}$) and $\hat{\alpha}(U)$ (for $KS_{jU}$).*

Next, we studied resampling for the three estimators. Subsampling with subsample size $m = 20 = 0.4n$ and the bootstrap was studied. Note that subsampling is consistent for all three estimators, but the bootstrap is consistent only for $AMA$. Panels (a) ((b)) of Figure 2, respectively, present the Kolmogorov–Smirnov distance, scaled by 100, between the distributions of $n^{1/2}(\hat{\alpha}_j - \alpha)$ and its subsampling (bootstrap) version, $j = MS, BMA, AMA$. We present the graphs for $|\beta| \leq 0.4 \approx 3/\sqrt{n}$, since there is not much difference between the three graphs for other values of $\beta$. It can be seen that the distances between the actual distribution and its subsampling/bootstrap versions are much smaller for $AMA$, while the resampling approximations for $MS$ and $BMA$ are particularly bad in the regions $\{|\beta| \in (0.1, 0.3)\}$. Also, there is little visual difference between the accuracies of the subsampling and the bootstrap approximations despite their different asymptotic behavior, which confirms some of the observations made in [1], [2] and [18].

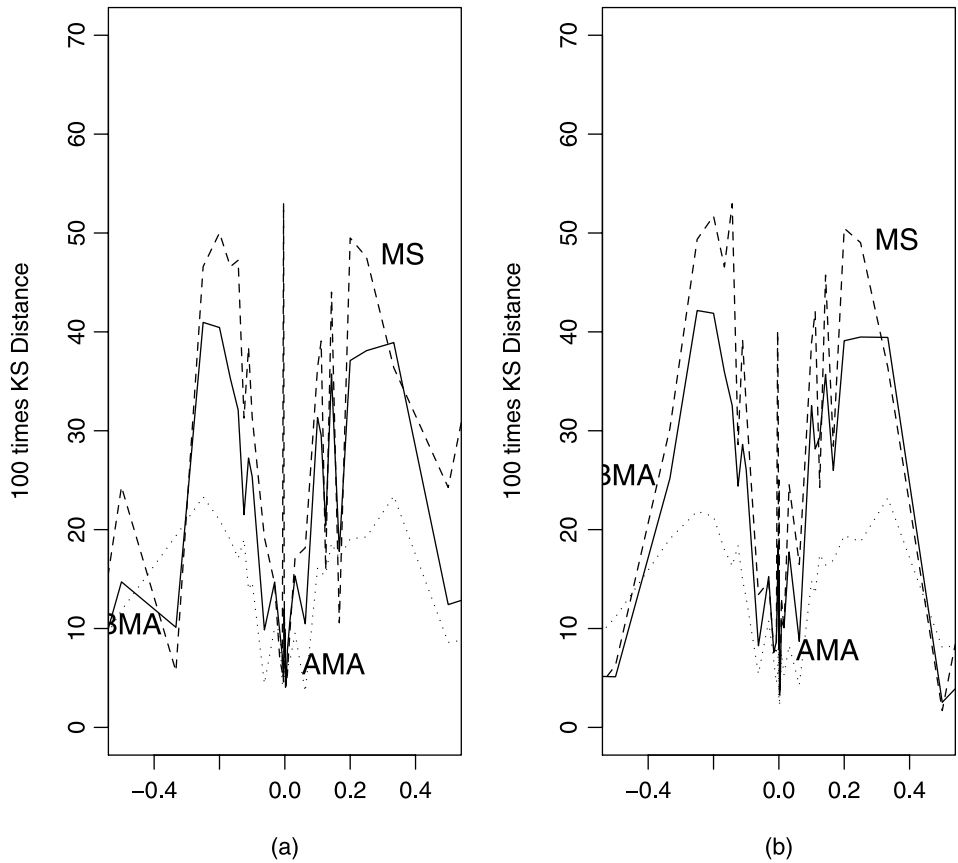

FIG 2. *Panel* (a) *is the subsampling approximation (subsample size* 20*) for the distribution of centered and scaled* $\hat{\alpha}_{BMA}$ *(solid line),* $\hat{\alpha}_{MS}$ *(broken line),* $\hat{\alpha}_{AMA}$ *(dotted line). Panel* (b) *is the corresponding bootstrap approximation.*

## 6. Discussion and conclusions

The problems associated with post-model-selection estimation have been discussed by several researchers. In current statistical practice, the process of selecting a model has similarities with hypothesis testing. On the other hand, estimation of parameters, some of which may be known constants in some of the models, is generally entirely separated from model selection. Estimation and testing/selection are two different paradigms of statistical analysis that are hard to integrate. The lack of uniformity across models that parameter estimators generally display, and the issues that arise subsequently, are products of the less than successful attempt to combine the two processes of estimation and selection.

In the Bayesian paradigm, model averaging seems to be a good integration of the two, since the selection step here is also an estimation exercise in spirit. The statement about integrated risks in Proposition 3.1 implies that Bayes' risks of model-averaged estimators are bounded. Thus, while minimaxity seems to be an elusive goal under model uncertainty, a fully Bayesian approach to analyzing risk behavior may be more successful.

In the context of bootstrapping model-averaged estimators, an alternative to $\check{\alpha}$

is to estimate the bias in $\hat{\alpha}$ in all the models, and define a bias corrected average of these. As the bias of $\hat{\alpha}(R)$ is $\beta ||X_1||^{-1} <X_1, X_2>$, if we estimate this by $\hat{\beta}||X_1||^{-1} <X_1, X_2>$, we get back $\hat{\alpha}(U)$. Nevertheless, in more complex problems the "bias corrected model averaged" estimator may be an interesting object to study.

In Theorem 4.1 we established the consistency of the paired bootstrap for a data-adaptive model-averaged estimator. Two other kinds of bootstrap are available in the linear regression context; namely, parametric bootstrap and the residual-based bootstrap. When only one model is in use, the parametric bootstrap generates data from it using estimated values for the unknown parameters, while the residual bootstrap obtains residuals after fitting the model. The equivalents of these are not obvious under model uncertainty.

In Section 4 we remarked that the data adaptive weights $p_{nR}$ and $p_{nU}$ may not share the same properties as the posterior model probabilities $\pi_{nR}$ and $\pi_{nU}$ of Section 3. It would be interesting to study when $p_{nR}$ and $p_{nU}$ can be interpreted as posterior probabilities, and also under what conditions the frequentist properties of a Bayesian model-averaged estimator may be elicited using the bootstrap.

**Acknowledgments.** Professor Chatterjee's research was partially supported by a grant from the University of Minnesota. We thank the referee and the editors of this monograph for some excellent comments and suggestions. Also, we would like to thank Professor Yuhong Yang, who carefully read an earlier draft of this paper and made several comments; which, along with several illuminating discussions, greatly enhanced our understanding on the scope and issues relating to model selection/averaging.